%% file: arxiv_r0.tex
\pdfoutput=1
\documentclass[11pt]{article}
\usepackage[margin=1in]{geometry}
\usepackage[utf8]{inputenc}
\usepackage[T1]{fontenc}
\usepackage{amsmath,amssymb,amsthm}
\usepackage{authblk}
\usepackage{orcidlink}
\usepackage[
  backend=biber,
  style=numeric-comp,
  sorting=none,
  maxbibnames=99
]{biblatex}

\newcommand{\minimize}{\operatornamewithlimits{minimize}}
\newcommand{\gate}{\kappa}

\usepackage{xcolor,doi}
\usepackage{hyperref}
\usepackage{listings}
\usepackage{booktabs}
\usepackage{tikz}
\usepackage{cleveref}
\usetikzlibrary{arrows.meta,backgrounds,calc,fit,matrix,positioning,shapes.geometric}
\input{figures/diagram_styles.tex}

\newtheorem{theorem}{Theorem}[section]
\newtheorem{proposition}[theorem]{Proposition}

\theoremstyle{definition}
\newtheorem{example}[theorem]{Example}
\hypersetup{colorlinks=true, linkcolor=blue!80!black, urlcolor=blue!80!black, citecolor=blue!80!black}
\usepackage[colorinlistoftodos]{todonotes}

\lstdefinestyle{pythonapi}{
  language=Python,
  basicstyle=\ttfamily\footnotesize,
  keywordstyle=\bfseries\color{blue!65!black},
  commentstyle=\itshape\color{green!35!black},
  stringstyle=\color{orange!60!black},
  emph={record,replay,grad,value_and_grad,replay_grad,replay_value_and_grad,replay_value_and_grad_batch,all_value_and_grad},
  emphstyle=\color{purple!70!black},
  emph={[2]jax,jnp,branch_extended_AD,bead,bnp},
  emphstyle={[2]\color{teal!70!black}},
  emph={[3]atol,rtol,tol_mode,input_scaled,abs_policy,paths,path},
  emphstyle={[3]\color{red!65!black}},
  numbers=left,
  numberstyle=\scriptsize\color{black!45},
  numbersep=0.75em,
  xleftmargin=2.5em,
  frame=none,
  backgroundcolor=\color{black!3},
  showstringspaces=false,
  columns=fullflexible,
  keepspaces=true,
  breaklines=true,
  breakatwhitespace=true,
  tabsize=2
}

\begin{document}

\title{Branch-Extended Automatic Differentiation}

\author[1]{Jan H\"uckelheim\,\orcidlink{0000-0003-3479-6361}\thanks{jhuckelheim@anl.gov}}
\author[1]{Jeffrey Larson\,\orcidlink{0000-0001-9924-2082}\thanks{jmlarson@anl.gov}}
\author[1]{Matt Menickelly\,\orcidlink{0000-0002-2023-0837}\thanks{mmenickelly@anl.gov}}
\affil[1]{Argonne National Laboratory, Lemont, Illinois, USA}
\date{}

\maketitle

\begin{abstract}

We present branch-extended automatic differentiation (BEAD), a software package that provides information about inactive branches of piecewise-differentiable functions, intended for use in nonsmooth optimization algorithms. BEAD is implemented in Python, building on the popular JAX framework. For supported nonsmooth operations such as abs and max, BEAD identifies alternative branches within tolerance from the given input, and computes gradients of the resulting branches extended to the given input. Such information has been previously used for robust optimization methods, but has required error-prone and tedious manual differentiation, unlike BEAD, which is fully automated for functions written in NumPy-like JAX code.
Our experiments show that our approach is as accurate as manually written derivative code, and offers more flexibility in modifying the function and choosing between branch-discovery rules. While we were originally motivated by the manifold sampling algorithm for nonsmooth composite objectives, we also integrate BEAD branch-extension gradients into a nonsmooth BFGS--SQP solver. On two small nonsmooth test problems, the resulting path-aware search resolves failures or substantially reduces the number of evaluations required by the BFGS--SQP solver.

\end{abstract}


\section{Introduction}
Automatic differentiation (AD) applies the chain rule of calculus to functions expressed as
computer programs and is a standard way to obtain derivatives without manually writing
or maintaining separate derivative code \cite{GriewankWalther2008}. 
For a
smooth program, differentiating its execution trace produces the derivative of
the represented function (with some caveats in the presence of roundoff, iterative solvers, and other implementation details, see~\cite{https://doi.org/10.1002/widm.1555} for an overview). 
Scientific codes, however,
often involve continuous but nonsmooth operations such as an absolute value,
a pointwise minimum or a pointwise maximum. 
When such operations appear as primitives in a function, the function typically becomes piecewise smooth; 
derivatives do not exist in the mathematical sense at inputs that lie exactly on a kink in the piecewise smooth structure.
In a computer program, such nonexistence is often represented at a low abstraction level as a \emph{branch} that switches control flow at such inputs.
At higher abstraction levels, this nonexistence occurs at \emph{branch points}, an intermediate (set of) value(s) at which a primitive enters one of several mathematically valid branches in defining the primal function.

As a motivating example of branches and branch points, consider an implementation of the maximum of two scalars $a$ and $b$:
\begin{verbatim}
    def max(a, b):
        if(a>=b) return a
        else return b
\end{verbatim}
This code consistently returns $a$ when $a=b$, an arbitrary branch: an equally
valid implementation returns $b$ instead:
\begin{verbatim}
    def max(a, b):
        if(a<=b) return b
        else return a
\end{verbatim}
For every instance where the \texttt{max} primitive appears in a program, the forward execution consistently follows one of two branches.
Conventional AD evaluates the derivative of the primal function determined by the corresponding execution trace.
However, users of the derivative program, especially practitioners of optimization, may want information about alternative branches not taken by the implementation.
We will clarify our intentions later, but note that the execution trace associated with alternative branches may have vastly different derivatives
\cite{LeeEtAl2020,BoltePauwels2020}.
A user may also want derivatives corresponding to alternative branches that were nearly followed within a prescribed tolerance; that is, for small $\texttt{eps}>0$,
replacing the first implementation of the \texttt{max} primitive with
\begin{verbatim}
    def max(a, b):
        if(a + eps >=b) return a
        else return b
\end{verbatim}
would have admitted an alternative branch for the same function input.

To overcome this limitation of conventional AD, we introduce \emph{branch-extended automatic differentiation} (BEAD).  At a given function input, BEAD records the branch taken at supported nonsmooth primitives such as \texttt{abs} and \texttt{max},
identifies alternative branches within a prescribed tolerance, and replays
the program once per \emph{path} through branches---an assignment of one branch at every branch point---thereby locking those branches.
Each replay defines a smooth
\emph{branch extension} at the query input, differentiable by
ordinary JAX AD~\cite{jax2018}.
The output is a set of \emph{branch-extension gradients} rather than only the gradient of the default path.
BEAD returns tolerance- and program-dependent branch-extension gradients;
these need not form a classical generalized derivative, as discussed in
\Cref{sec:generalized-derivatives}.

Our implementation and test cases are mostly motivated by manifold sampling, a derivative-free optimization method that relies on composite objective functions of a certain form, discussed further in \Cref{sec:background}.
However, BEAD's applicability is not specific to composite objectives or derivative-free optimization, and we demonstrate some preliminary results with scalar-valued nonsmooth functions in the context of the constrained
nonsmooth solver PyGRANSO \cite{liang2022ncvx,curtis2017bfgssqp} in \Cref{sec:pygranso-cases}.
BEAD supports both local tolerances on intermediate values and input-scaled
tolerances based on first-order input sensitivities, as described in
\Cref{sec.discovery}.

The contributions of this work are a record-and-replay formulation of BEAD, a
JAX implementation for array operations that commonly appear in definitions of piecewise smooth functions, implementations of two branch-discovery rules, and an evaluation
against manual implementations of branch-extension gradients from past work.

\section{Background}
\label{sec:background}
We now review piecewise-smooth programs and existing generalized-derivative
concepts. 
We then review BEAD's connection to nonsmooth optimization methods.

\subsection{Piecewise-smooth programs}

A continuous piecewise-$C^1$ function is locally represented by finitely many
smooth selection functions, with one or more selection functions active at a switching
surface \cite{Scholtes2012}. 
For input queries off the switching surfaces, ordinary AD of the
execution trace returns the derivative of the locally unique selection function. 
For input queries on the switching surface, a single execution trace corresponds to only one selection function; it does not describe the other selection functions that determine nearby
first-order behaviors of the piecewise-smooth function.

In related work, abs-normal form makes this switching structure explicit by representing
piecewise-smooth functions using smooth operations connected by absolute-value
switches \cite{Griewank2013AbsNormal}.  It supports systematic
piecewise-linearization and generalized algorithmic differentiation.  BEAD is
closely related but does not transform a function into abs-normal
form, nor does it return a complete piecewise-linear local model of the function. Instead, it records the branch taken at each supported JAX primitive, identifies within-tolerance
alternative branches, extends them to the current input, and applies ordinary AD to each, producing branch-extension values and gradients.

We list several other related works in piecewise-smooth programs here. 
Intensional derivatives formalize
the derivatives computed by AD for piecewise-analytic programs and agree with
ordinary derivatives almost everywhere \cite{LeeEtAl2020}.  Conservative
fields provide a chain-rule-compatible set-valued calculus for broad classes of
nonsmooth programs and support convergence analyses for stochastic-gradient
methods \cite{BoltePauwels2021,BoltePauwels2020}.  
Automatic
subdifferentiation seeks to compute a mathematically valid generalized
subderivative for a qualified program \cite{KakadeLee2018}.  
These approaches
primarily justify and select generalized derivative information, but do not provide information about alternative branches, which is required for robust nonsmooth optimization solvers that motivate this work.

\subsection{Relationship to generalized derivatives}\label{sec:generalized-derivatives}
Several notions of generalized derivative exist~\cite{Clarke1990,RockafellarWets1998}.
These generalized derivatives include Dini or Hadamard directional
derivatives that describe a one-sided rate of change in a specified direction.
The
Bouligand, or B-, subdifferential $\partial_B f(x)$ collects limits of Jacobians from
differentiable points in a neighborhood of a given point.
For locally Lipschitz
scalar functions, the Clarke subdifferential $\partial_C f(x)$ and Clarke generalized directional derivative represent the convex hull of limiting
gradients and their support function, respectively.
Convex subgradients are defined by global supporting inequalities for the epigraph of a convex function. 

In general, BEAD computes none of these generalized derivatives.
Given an input query $x$, BEAD produces a finite set
\[
\mathcal G_\tau(x)=\{\nabla f_p(x):p\in\mathcal P_\tau(x)\},
\]
where $\tau>0$ is an absolute tolerance used in determining alternate branches,
$\mathcal P_\tau(x)$ is the set of recorded paths through branches, and
$f_p$ is the branch extension obtained by locking path $p$.
The set $\mathcal P_{\tau}(x)$ depends on
the program, supported primitives, and tolerance strategy. 
\Cref{sec:tau-zero} gives theoretical conditions under which
$\mathcal G_0(x)$ contains $\partial_B f(x)$.

Moreover, with a positive tolerance $\tau>0$, inactive smooth extensions are deliberately
included in $\mathcal P_\tau(x)$, and so it is generally not the case when
$\tau>0$ that $\mathcal G_{\tau}(x)$ coincides with any of these generalized
derivatives.

Another related approach is lexicographic differentiation~\cite{Nesterov2005Lexicographic}, which computes directional derivatives for a provided sequence
of infinitesimal perturbations.
Lexicographic differentation
has been employed to find solutions of both systems of nonsmooth equations and parametric differential
equations \cite{KhanBarton2015GeneralizedDerivatives}.
The directional nature of lexicographic differentiation distinguishes it from BEAD, which instead discovers a set of 
paths according to a finite tolerance and differentiates each corresponding branch extension separately.

\subsection{Relationship to nonsmooth optimization methods}\label{sec:nonsmooth-optimization}
Various methods have been proposed for the solution of nonsmooth optimization problems, that is, optimization problems for which the objective function and possibly some constraints are nonsmooth functions of their inputs. 
The most robust of these methods generally employ a generalized derivative at points of nondifferentiability, or else construct approximations of generalized derivatives from conventional gradients at slightly perturbed points in neighborhoods of points of nondifferentiability. 
Bundle methods~\cite{Lemarechal1975,Wolfe1975,HiriartUrrutyLemarechal1993} use subgradients at various query points in the optimization domain to build and update cutting-plane models of objectives and constraints. 
Gradient sampling methods~\cite{BurkeLewisOverton2005,burke2020gradient} evaluate
conventional gradients at randomly sampled points of differentiability near a given input and
computes the gradients' convex hull to approximate a Clarke subdifferential. 

In contrast to the mechanisms of bundle methods and gradient sampling methods, both of which aggregate (sub)gradient information across multiple input points, BEAD considers a single input point and
differentiates branch extensions identified from program structure.
This mode of computation is particularly useful within the paradigm of manifold
sampling~\cite{Larson2016c,KLW2018,Larson2020,Larson2022}, which seeks to identify (a subset of) active manifolds defining a continuous selection at any given incumbent. 
Our motivating manifold-sampling problems have the composite form
\begin{equation} \label{eq:prob_def}
\minimize_{x \in \mathbb{R}^n} f(x) = h(F(x)),
\end{equation}
where $F : \mathbb{R}^n \to \mathbb{R}^p$ is expensive to evaluate and
$h : \mathbb{R}^p \to \mathbb{R}$ is a known, possibly nonsmooth scalar
mapping.
Manifold sampling builds separate models of the components of $F$
and uses values and gradients associated with active manifolds of $h$;
``activity'' here generally corresponds to branch extensions.
BEAD automates
the path-identification and differentiation interface required by manifold sampling
methods. We discuss benchmark problems further in \Cref{sec:test-cases}.

However, BEAD can be applied beyond manifold sampling methods.
As a particular example, a finite BEAD set $\mathcal G_\tau(x)$ can be employed to construct a cutting-plane model that would appear in a bundle method. 
The downstream bundle method may minimize the support function
$\max_{g\in\mathcal G_\tau(x)} g^T d$ when choosing a direction for use in a line-search or in a trust-region subproblem. 
This use case does not require a decomposition
$f=h\circ F$ like manifold sampling is designed to exploit; instead, BEAD can be applied directly to a piecewise-smooth function and generate BEAD sets of practical utility.

\section{Implementation}
We now BEAD's record-and-replay execution model, its branch-discovery rules,
batched and compiled replay, and its user-facing API.

\subsection{Design and execution model}\label{sec:design}

BEAD is implemented as a Python package built on JAX~\cite{jax2018}.  
BEAD does not replace or modify JAX's differentiation machinery in any significant way.
Instead, BEAD changes the behavior of the forward pass, generating multiple
branch extensions that can then be differentiated normally.
User code defining the primal function must import a NumPy-like wrapper module for a set of supported nonsmooth
operations, as well as \texttt{jax.numpy} in order to define smooth operations.

A \emph{recording pass} executes the user's primal function once.
At each supported nonsmooth primitive, BEAD records the executed branch and identifies alternative
branches that are within a specified tolerance, to be defined in \Cref{sec.discovery}.
The implementation of BEAD distinguishes several types of primitives.
\emph{Reduction extrema} primitives, such as \texttt{min} and \texttt{max}, produce a scalar output and record the selected index together with other indices within tolerance of the extreme. 
A \emph{unary elementwise} primitive, such as \texttt{abs}, records a sign branch at entries within tolerance of zero.
\emph{Binary elementwise extrema} primitives, such as \texttt{minimum} and \texttt{maximum}, record an operand branch at each array index.
These branches are stored in small typed records and can have large downstream effects: elementwise primitives can induce exponentially many branch combinations when many entries are within tolerance simultaneously, whereas reduction extrema primitives produce a number of branches that grows, at most, linearly in the size of the input vector.
These combinatorial issues should be taken into account when designing primal functions.

The \emph{recorded branches}, that is the executed branches and alternative branches within tolerance taken together, define a \texttt{PathSet}.  
A single resolved path through branches within a \texttt{PathSet} fixes one
branch at every branch point in the program; downstream code can refer to it by a path label string.
The number of paths in a \texttt{PathSet} is generally exponential in the number of branch points that admit within-tolerance alternatives for a given input.
To avoid instantiating an exponential number of paths in memory whenever possible,
our implementation uniquely identifies each path by a single integer, which represents the position of that path in a sequence defined by the Cartesian product of the branches available at each branch point encountered in a program.
A \texttt{PathSet} by default therefore only needs to store the available branches at each branch point separately, as opposed to storing all possible combinations.

\texttt{PathSet} objects also support set operations such as \texttt{union} and \texttt{intersection}.
Such operations may destroy the Cartesian-product structure, in which case BEAD stores an explicit list of paths.  Whether stored compactly or explicitly, evaluating every path requires work proportional to the number of paths, which is exponential in the number of branch points that admit alternatives within tolerance.
Thus, users of BEAD may need to control the runtime or memory use of this approach either by judiciously reducing the tolerance or designing their problem in a way that ensures that the number of alternatives is sufficiently bounded. 
Future work might investigate, for example, automatic sampling strategies in BEAD in order to limit resource use when user-controlled remedies are apparently insufficient.

A \emph{replay pass} executes the user program using the branches fixed by
one resolved path in a \texttt{PathSet}.
Thus, a replay pass changes the forward values produced by the
supported nonsmooth primitives. 
Once those branches are locked, the replay
pass is an ordinary smooth JAX program for that branch extension, wherever its
elementary operations are smooth. 
BEAD can consequently apply
\texttt{jax.grad} or \texttt{jax.value\_and\_grad} directly to the replayed
function. 
In particular, BEAD does not implement a separate reverse pass or a
custom derivative rule for each composite nonsmooth function; by changing the
forward branches, JAX may compute the derivative of that
branch extension automatically. 

\begin{figure}[h!]
  \centering
  \input{figures/record_replay_design.tex}
  \caption{Record--replay example.  Recording identifies branches within tolerance
    at a maximum and an absolute value, producing four paths.  Replaying the
    highlighted path replaces the original forward computation by
    $f_{p_4}=g_2-h$, after which ordinary JAX automatic differentiation returns
    the value and gradient of that branch extension.  Multiple resolved paths can be evaluated by the compiled pipeline.}
  \label{fig:record-replay-design}
\end{figure}
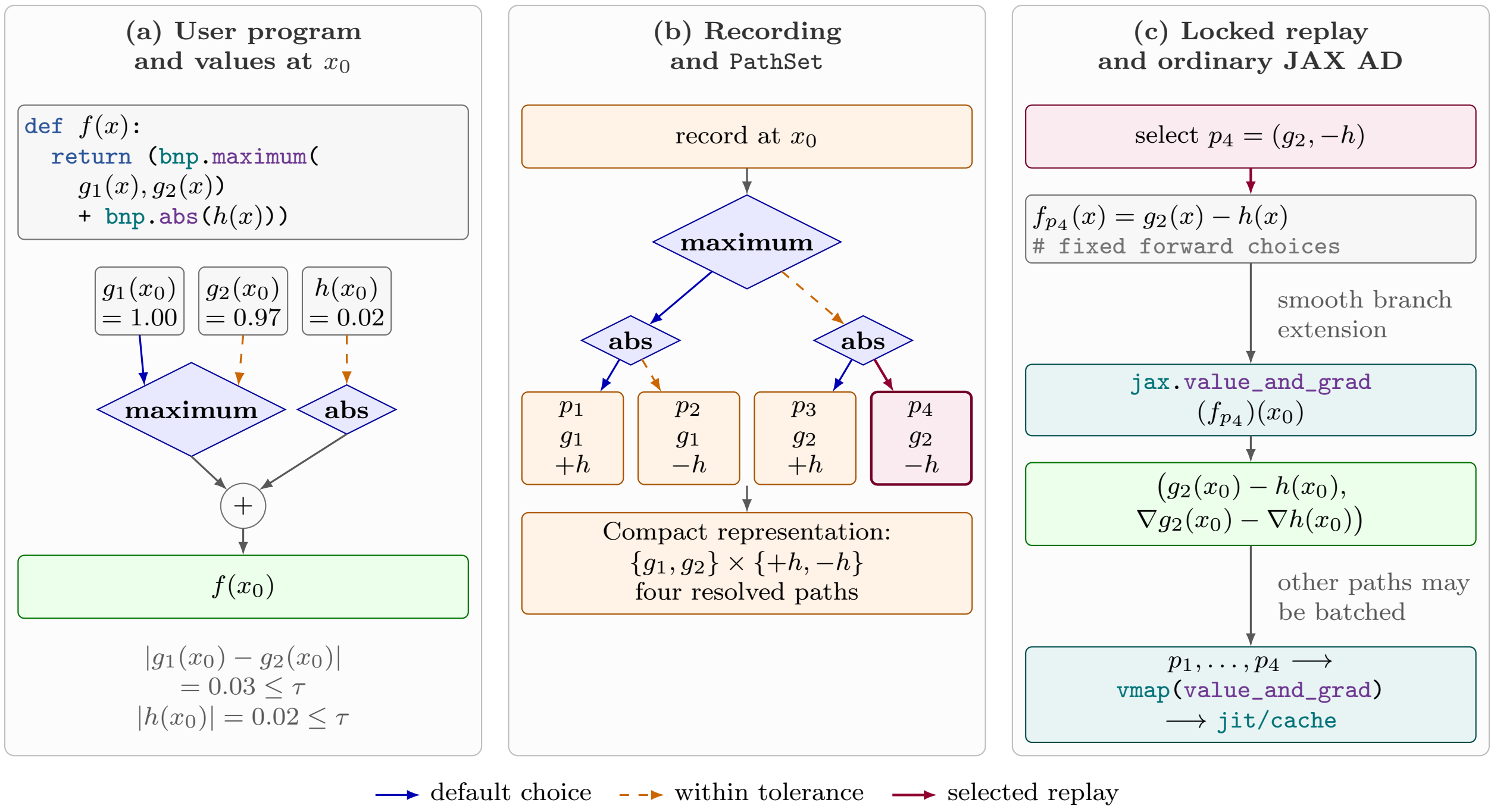

Figure~\ref{fig:record-replay-design} summarizes this execution model.  
Branch discovery is performed during forward execution, while each replay is an
ordinary JAX computation that can be differentiated and compiled.

\subsection{Branch discovery and path-set growth}\label{sec.discovery}
Branch discovery is controlled by absolute and relative tolerances.
We propose two explicit rules for branch discovery employing these tolerances.
Every supported primitive has an associated scalar switching condition $q$, computed
from the intermediate value(s) that primitive takes. 
The sign of $q$ determines the branch taken at that branch point.
For example, \texttt{abs} switches on the sign of its single input, while
\texttt{maximum} switches on the sign of the difference between its two
operands.
BEAD applies each primitive's $q$ consistently during both recording and
replay of branches.
Because each branch point is tested independently, a combination of branches
recorded in $\mathcal P_\tau(x)$ need not correspond to any single, jointly
reachable input perturbation; the path set could include branches that can't be
reached by any nearby input; \Cref{ex:tau-zero} is a small instance.
The local branch-discovery rule
admits the alternative branch when $|q|$ is below a tolerance expressed in
the units of that intermediate quantity.
We implement such a rule as
$$|q| \leq \mathtt{atol}+\mathtt{rtol}|v|,$$
where $v$ is an operation-specific scale: the selected extremum for a
reduction, the larger-magnitude operand for a binary extremum, and the input to
an absolute value.

Such a test is generally inexpensive and easy to implement, but may not always behave in a way that is intuitive to a user; intermediate values and their sensitivity to function input perturbations may have vastly different scales from one another.
As an alternative, it may be more intuitive to ask whether perturbing the
function input within some small tolerance could flip the sign of $q$.
In other words, we search for alternative program paths that would be taken for some input within a small neighborhood around the current function input $x$.  
Exact propagation of such a neighborhood through a
nonlinear program, however, is a challenging problem. 
Interval methods can propagate bounds through a program, but may overestimate reachable ranges.
Taylor models retain higher-order information at greater cost
\cite{MooreKearfottCloud2009,MakinoBerz2003}.

Noting these challenges, BEAD employs an \emph{input-scaled} branch-discovery rule to serve as a first-order approximation to the propagation of a neighborhood through a nonlinear program.
Because $q$ is computed from the function input $x$ through the forward pass,
it can be viewed as a scalar function $q(x)$ of $x$, with
switching surface $\{x:q(x)=0\}$.
For an input ball $\|\delta\|\leq\tau$, the first-order range of the switching
quantity is controlled by $\tau\|\nabla q(x)\|_*$, where $\|\cdot\|_*$ is the
dual norm.  Forward-mode AD propagates input sensitivities through the program;
BEAD uses the resulting magnitudes to estimate this scale and adjust the local
tolerance.
The sensitivities only guide branch discovery; replay gradients are still computed
by ordinary JAX AD. 
This sensitivity-based approximation does not provide any strong bounds, as it deliberately neglects higher-order effects for nonlinear functions. 
Like the local rule, the input-scaled rule is a practical branch-discovery rule and generally does not correspond precisely to a generalized derivative.

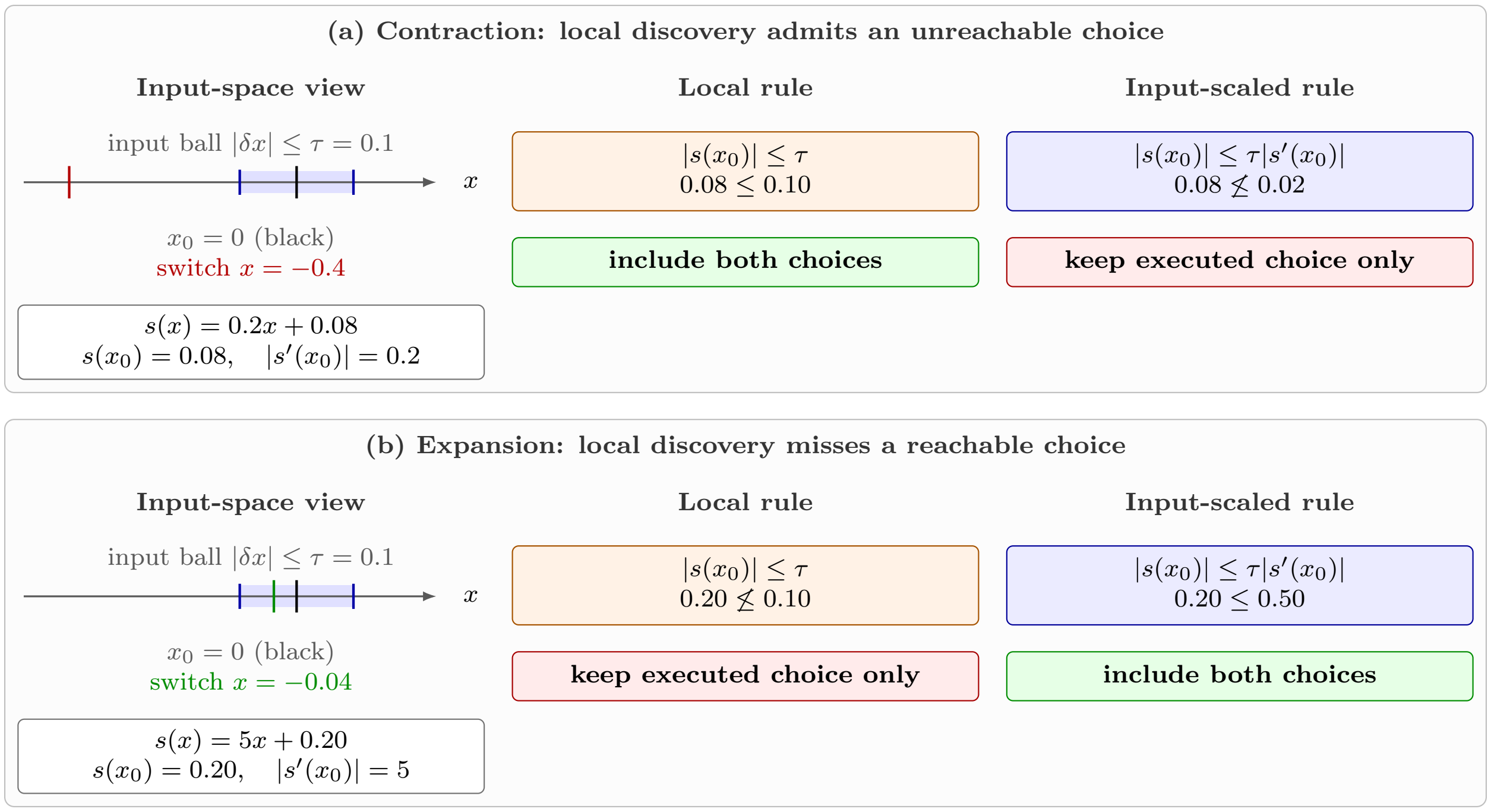
\begin{figure}[h!]
  \centering
  \input{figures/local_input_scaled.tex}
  \caption{Local and input-scaled branch discovery for two scalar linear
    switching quantities, using \texttt{atol}$=\tau=0.1$ and
    \texttt{rtol}$=0$.  In panel (a), contraction makes the local rule include
    an alternative whose switching point is outside the input ball.  In panel
    (b), expansion makes the local rule miss an alternative whose switching
    point is inside the input ball.  Scaling by $|q'(x_0)|$ gives the correct
    first-order decision in both examples because $x$ is scalar and the input
    norm is absolute value.}
  \label{fig:local-input-scaled}
\end{figure}

Nonetheless, the first-order input-scaled rule corrects two failure modes of
the local rule in the linear examples below.  More generally, if
$q$ has an $L$-Lipschitz gradient on the input ball, then
\[
 q(x+\delta)=q(x)+\nabla q(x)^T\delta+R(\delta),\qquad
 |R(\delta)|\leq \tfrac{L}{2}\|\delta\|^2.
\]
The input-scaled rule models the linear term and neglects the remainder.
If the switching quantity $q(x)$ \emph{contracts} input perturbations, the local rule may admit a
branch whose switching surface is outside the input ball.
If $q(x)$ \emph{expands} input
perturbations, the local rule may reject a branch whose switching surface lies
inside the ball.
The input-scaled rule corrects both decisions for these linear
examples.
We illustrate an example of these behaviors in~\Cref{fig:local-input-scaled}. 

As an example of the mitigation strategies mentioned above, we present a policy specific
to the \texttt{abs} function applied entrywise to a vector.
When the vector contains $k$ entries that are near zero within tolerance, there are $2^k$ recorded paths, enumerated by a $k$-fold Cartesian product of sign branches $\pm$.
Naively, one could employ a policy, called via
\texttt{abs\_policy="enumerate"} in BEAD, that causes BEAD to record both signs $\pm$ for each
ambiguous entry.
However, with the default \texttt{abs\_policy="zero"}, BEAD instead retains one
compact representative branch whose recorded value is simply the absolute value, and the associated derivative
is taken at zero at any ambiguous entry within tolerance $\tau$.
This zero is not one of the two one-sided derivatives $\pm 1$ of \texttt{abs};
rather, it is the minimum-norm element of its Clarke subdifferential $[-1,1]$
at a tie, the standard representative when a single (sub)gradient is required,
for example in a stationarity measure.
The default policy avoids $2^k$ recorded paths in favor of one path, at the expense of a less informative $\mathcal{G}_{\tau}(x)$.

The absolute-value policy separates three decisions: which alternatives are
discovered, how they are represented, and which resolved paths are evaluated.
Other primitives with many within-tolerance alternatives could make the same
separation; discovery, representation, and evaluation need not all be
exhaustive.
Future work involves topics such as lazy path generation, sampling,
limits on retained paths, removal of duplicate or nearly identical
gradients, compact representations of independent or convex derivative sets,
and optimizer-guided selection of paths relevant to a search direction.
Such policies should be studied in terms of downstream accuracy (e.g., in
stationarity tests or optimization search directions), memory use, and runtime.
Currently, \texttt{abs} is the only primitive that has such a policy implemented in BEAD, because other primitives did not demonstrate path-set explosion in our tests. Future work will study similar policies, for example, for the \texttt{maximum} and \texttt{minimum} primitive.
These primitives in particular exhibit an explosion in the size of discovered path-sets when two compared input vectors have many indices that are equal in value or nearly equal in value.

\subsection{Generalized derivatives at zero tolerance}\label{sec:tau-zero}

We consider the special case of $\mathcal{G}_{\tau}(x)$ where $\tau=0$.
Every recorded path in $\mathcal{P}_0(x)$ fixes one branch at each 
branch point, and so defines a smooth branch extension $f_p$ that agrees with $f$ at $x$.
These branch extensions are selection functions in the sense of Scholtes~\cite{Scholtes2012}.
To formalize a result about $\mathcal{G}_0(x)$, we call a path (selection function) $p\in\mathcal P_0(x)$ \emph{essentially
active} if $x$ lies in the closure of the interior of $\{y: f(y)=f_p(y)\}$.

\begin{proposition}\label{prop:tau-zero}
Let every operation in the program defining $f$ be either continuously
differentiable or else a supported nonsmooth primitive. 
Suppose branch discovery with $\tau=0$ successfully records
all alternative branches at every branch point. 
That is, branch discovery records every branch that the program can execute at inputs arbitrarily close to $x$.
Suppose further that
$\mathcal P_0(x)$ is the Cartesian product of the recorded branches, every
path of which is evaluated.  
Then
\begin{equation}\label{eq:inclusion}
  \partial_B f(x)\subseteq\mathcal G_0(x)
  \quad\text{and}\quad
  \partial_C f(x)\subseteq\operatorname{conv}\mathcal G_0(x),
\end{equation}
and equality holds in the first inclusion if and only if every
$p\in\mathcal P_0(x)$ satisfies $\nabla f_p(x)=\nabla f_{p'}(x)$ for some
essentially active $p'\in\mathcal P_0(x)$.
In that case, equality holds in the second inclusion as well.
\end{proposition}

\begin{proof}
A recorded path fixes a branch at every branch point, so each $f_p$ is a
composition of continuously differentiable operations and is $C^1$ in some open neighborhood of $x$.
Let $U_1$ be the intersection of these finitely many neighborhoods.

We will show that $f$ is a selection of the family $\{f_p\}$ near $x$.
Consider the execution of the program defining $f$. At any
branch point whose switching quantity $q$ is nonzero at the intermediate values 
entering the branch point (call these values $v(x)$), the same branch is executed for all values in an open neighborhood of $v(x)$. 
By continuity of $f$, this means the same branch will be executed for all $y$ in some open neighborhood of $x$. 
Denote the intersection of these finitely many (corresponding to each branch point with nonzero intermediate values) open sets by $U_2$. 
On the other hand, at any branch point whose switching quantity $q$ vanishes at $v(x)$, every alternative is recorded by assumption.
Consider the open set $U := U_1\cap U_2$. 
We have both that $x\in U$ and
for every $y\in U$, the branches executed by the program at $y$ belong to some path $p\in\mathcal P_0(x)$, that is,  $f(y)=f_p(y)$.
Thus $f$ is a continuous selection of the finite family of $C^1$ functions $\{f_p\}_{p\in\mathcal P_0(x)}$ on $U$.

\cite[Prop.~4.3.1]{Scholtes2012} gives that $\partial_B f(x)$ is the set of
gradients at $x$ of the essentially active members of this family, and that
$\partial_C f(x)$ is its convex hull.
Because every path in the full product is evaluated,
$\mathcal G_0(x)=\{\nabla f_p(x):p\in\mathcal P_0(x)\}$, which contains the gradients of the essentially active paths. This is the first inclusion in \eqref{eq:inclusion}; taking convex hulls gives the second.
Finally, $\mathcal G_0(x)=\partial_B f(x)$ exactly when every $\nabla f_p(x)$ equals the gradient of some essentially active path, which is the stated criterion. Equality of the sets implies equality of their convex hulls.
\end{proof}

The following simple one-dimensional example shows the hypotheses of \Cref{prop:tau-zero} can hold while equality in \cref{eq:inclusion} fails.

\begin{example}\label{ex:tau-zero}
For $f(t)=t+\max\{t,0\}-\max\{t,0\}\equiv t$ the two
maxima form independent branch points in the program.
Differentiating the four resulting paths gives
$\mathcal G_0(0)=\{0,1,2\}$ against $\partial_B f(0)=\{1\}$.  Here
$0\in\operatorname{conv}\mathcal G_0(0)$, so a method measuring stationarity by
$\operatorname{dist}(0,\operatorname{conv}\mathcal G_0(x))$ would report
$x=0$ as stationary.
\end{example}

The hypotheses of \Cref{prop:tau-zero} hold in BEAD when \texttt{atol}$=$\texttt{rtol}$=0$,
\texttt{tol\_mode="local"}, and \texttt{abs\_policy="enumerate"}.%
With the default \texttt{abs\_policy="zero"}, they fail:
$\mathcal G_0(0)=\{0\}$ for $f(t)=|t|$, against $\partial_B f(0)=\{-1,1\}$.
The guarantees of \Cref{prop:tau-zero} also fail for $\tau>0$: a recorded branch need no longer agree with
any executed branch at $x$, and so $f_p(x)\neq f(x)$ is possible; the extensions then no
longer form a selection family for $f$.

Note that BEAD does not test essential activity, which \Cref{prop:tau-zero} requires.
For programs built from smooth operations and
absolute values, the linear independence kink
qualification (LIKQ) of abs-normal form \cite[Def.~2]{GriewankWalther2016} is a sufficient condition: under it,
every sign combination at the active switches is realized on a set of positive
measure \cite[Cor.~7]{BaumgaertnerBethke2024}, so every recorded path is essentially active and equality holds in \eqref{eq:inclusion}.%

\subsection{Batched and compiled replay}

JAX supports JIT (just-in-time) compilation to accelerate the execution speed of functions and their gradient computation. 
JIT compilation introduces a tradeoff.  
The first call to JAX incurs tracing and
compilation costs, whereas later calls can reuse the compiled program and avoid
repeated Python dispatch. 
JAX usually needs to re-compile programs when their control flow or certain array dimensions change. 
To avoid frequent recompilations during the interactions of an optimization method with BEAD, our implementation pads the size of some data structures to the next power of two and removes padded results afterward. 
This performs some unused work but
reduces the number of compiled batch shapes.
Compiled replay functions are
cached by the user function, operation layout, differentiated arguments, and
auxiliary-output option.  
This optimization generally reduces replay overhead but is independent of the possibly combinatorial enumeration of recorded paths.

We call a collection of resolved paths that are evaluated together a \emph{path
batch}; a path batch is evaluated in \emph{batched replay}.
For a
collection of paths with the same sequence of branch points, BEAD first
validates their structure and converts their corresponding branches into dense JAX arrays.
BEAD
then applies \texttt{jax.value\_and\_grad} to one replay, uses \texttt{jax.vmap}
to evaluate the resolved paths together, and applies \texttt{jax.jit} to compile
the vectorized calculation.

\subsection{User-facing API}

The public interface follows JAX's functional-transformation style.
That is, each BEAD
transformation accepts a Python function and returns a callable with the same
user arguments. 
Supported nonsmooth operations are written using
\texttt{branch\_extended\_AD.numpy}; smooth operations may use
\texttt{jax.numpy}. 
Table~\ref{tab:bead-api} summarizes the layers of the API.

\begin{table}[htbp]
  \centering
  \caption{Main BEAD transformations.}
  \label{tab:bead-api}
  \begin{tabular}{ll}
    \toprule
    User goal & Transformation \\
    \midrule
    Discover within-tolerance paths & \texttt{record} \\
    Evaluate one resolved path & \texttt{replay} \\
    Differentiate the default path & \texttt{grad}, \texttt{value\_and\_grad} \\
    Differentiate one resolved path & \texttt{replay\_grad}, \texttt{replay\_value\_and\_grad} \\
    Differentiate all discovered paths & \texttt{all\_value\_and\_grad} \\
    Evaluate known paths in one compiled call & \texttt{replay\_value\_and\_grad\_batch} \\
    \bottomrule
  \end{tabular}
\end{table}

Listing~\ref{lst:api-basic} provides a demonstration of recording and replay. 
A recording pass returns
the value for the default path along with a \texttt{PathSet}.  
Replay accepts one resolved path and
returns the corresponding branch-extension value.

\begin{lstlisting}[style=pythonapi,caption={Recording and replaying within-tolerance paths.},label={lst:api-basic}]
import jax.numpy as jnp
import branch_extended_AD as bead
import branch_extended_AD.numpy as bnp

def h(z):
    return bnp.sum(bnp.maximum(z, 0.0))

z = jnp.array([0.0, 1.0, -2.0])
value, paths = bead.record(h, atol=1.0e-8)(z)

for path in paths:
    locked_value = bead.replay(h, path)(z)
\end{lstlisting}

Derivative transformations use the same path objects.  The default path
routines return the default-path derivative together with the paths discovered at
the provided input. 
Replay routines evaluate one known path.  When every recorded
path is required, \texttt{all\_value\_and\_grad} performs discovery and replay
in one call.

\begin{lstlisting}[style=pythonapi,caption={Computing branch-extension gradients from resolved paths.},label={lst:api-grad}]
(value, default_grad), paths = bead.value_and_grad(h, atol=1.0e-8)(z)
branch_data = [bead.replay_value_and_grad(h, path)(z) for path in paths]
all_branch_data, paths = bead.all_value_and_grad(h, atol=1.0e-8)(z)
\end{lstlisting}

The options \texttt{atol}, \texttt{rtol}, \texttt{tol\_mode}, and
\texttt{abs\_policy} dictate the scope of branch discovery. 
The options \texttt{argnums} and
\texttt{has\_aux} follow the corresponding JAX transformations.  
If a user employing BEAD has already collected resolved paths, they may use
\texttt{replay\_value\_and\_grad\_batch} to evaluate the same resolved paths through the compiled
batched replay rather than needing to rediscover them.

\section{Test cases}
\label{sec:test-cases}
We now discuss several use cases designed to probe the capabilities of BEAD, and note the capabilities' alignment with our design intentions. 

\subsection{Manifold sampling}

All manifold-sampling test cases use the composite form
\eqref{eq:prob_def}.  
A traditional derivative-free method will generally build and maintain one
surrogate $m_f$ to model the full objective $f = h\circ F$. 
Manifold sampling instead builds and maintains 
separate models $m_1,\ldots,m_p$, aggregated in the vector-valued model $m_F$, for the components $F_1,\dots, F_p$ of $F$ and combines these $p$ models
using the known structure of $h$. 
The componentwise construction $h\circ m_F$ often provides
a more accurate local model of $f$ than $m_f$. 

When $h$ is nonsmooth, a manifold sampling method requires both function value and gradients for its active
manifolds.
For a rigorous definition of active manifold, we refer the reader again to~~\cite{Larson2016c,KLW2018,Larson2020,Larson2022}. 
As a motivating example of the definition for the sake of understanding, for $h(z)=\lVert z\rVert_1$, manifold identification reduces to
labeling the signs of the entries of $F(x)$.
For a more complicated mapping,
such as the censored-$\ell_1$ loss
\[
 h(F(x))=\sum_{i=1}^{p}\left|d_i-\max\{F_i(x),c_i\}\right|,
\]
manual identification and differentiation are substantially more involved.
BEAD automates this problem-specific work and reduces the risk of errors.

To test whether our AD-encoded $h$ functions are correct, we needed 
a test suite broad enough to exercise BEAD's various tolerances and branch discovery policies. 
Rather than hand-selecting a handful of problem/$h$ combinations, we
ran every smooth vector-valued function $F$ from a benchmark set against every
$h$ we had previously hand-coded in the IBCDFO package~\cite{ibcdfo}.

For $F$, we use the full set of 53 test problems from the Mor{\'e}--Wild
DFO benchmark collection~\cite{JJMSMW09}; problems include a reference starting
point $x_0$. 
These problems span a low
range of dimensions ($2\le n \le 12)$ and output dimensions $2 \le p \le 45$.
The Mor{\`e}--Wild benchmark is intentionally small in scale, because it is intended to represent the types of problems typically encountered in simulation-based optimization, where $F$ is very expensive to evaluate and so a single simplex gradient approximation (requiring $n+1$ function evaluations) is cost-prohibitive. 

As for $h$, to really test the capabilities of BEAD, we defined nine different $h$ functions, listed in Table~\ref{tab:hfuns}.
For many of these $h$ functions (e.g., piecewise min/max functions) optimization
regularly encounters points of nondifferentiability, making branch discovery of critical importance.
Other $h$ functions (e.g., the piecewise-quadratic and censored-$\ell_1$
losses) are compositions with so many branching locations that
hand-differentiation can be impractical and error-prone.

\begin{table}[htbp]
\centering
\caption{The nine $h$ functions used in the benchmark, each applied to $z = F(x) \in \mathbb{R}^p$.}
\label{tab:hfuns}
\begin{tabular}{ll}
\toprule
Name & Formula \\
\midrule
\texttt{one\_norm} & $h(z) = \sum_{i=1}^p |z_i|$ \\
\texttt{pw\_maximum} & $h(z) = \max\{z_1, \ldots, z_p\}$ \\
\texttt{pw\_maximum\_sq} & $h(z) = \max\{z_1^2, \ldots, z_p^2\}$ \\
\texttt{pw\_minimum} & $h(z) = \min\{z_1, \ldots, z_p\}$ \\
\texttt{pw\_minimum\_sq} & $h(z) = \min\{z_1^2, \ldots, z_p^2\}$ \\
\texttt{quantile} & $q$-th quantile of $\{z_1^2, \ldots, z_p^2\}$ \\
\texttt{max\_plus\_qvp} & $\max\{z_1, \ldots, z_{p-1}\} + \alpha \sum_{i \ge p} \max\{z_i, 0\}^2$ \\
\texttt{censored\_L1\_loss} & $h(z; c, d) = \sum_{i=1}^p |d_i - \max(z_i, c_i)|$ \\
\texttt{piecewise\_quadratic} & $h(z; Q_j, z_j, c_j) = \max_j \left\{ \|z - z_j\|^2_{Q_j} + c_j \right\}$ \\
\bottomrule
\end{tabular}
\end{table}

Two of the $h$ functions in \Cref{tab:hfuns} deserve discussion. 
In our tests,
\texttt{h\_quantile} is hardcoded to $q=1$, which makes it functionally
identical to \texttt{h\_pw\_minimum\_squared}; we therefore exclude \texttt{h\_quantile} from
the counted instances below to avoid redundancies.
We keep \texttt{h\_quantile} in Table~\ref{tab:hfuns} only to
document the general quantile-based $h$ function. 

Similarly,
\texttt{h\_max\_plus\_quadratic\_violation\_penalty} has its penalty
weight $\alpha$ hardcoded to $0$, and so it reduces to
\texttt{h\_pw\_maximum} restricted to the first $p-1$ components; we choose to keep this $h(z)$ in the test set (and in the instance count) because it exercises a
different code path even though its values coincide with
\texttt{h\_pw\_maximum}.

The last two entries in Table~\ref{tab:hfuns} are not fixed functions but
function factories: \texttt{create\_\allowbreak censored\_L1\_loss\_hfun} builds a one-sided
censored $\ell_1$ loss from problem-specific data $(c, d)$, following
Womersley's censored discrete linear $\ell_1$ approximation
setup~\cite{Womersley_1986}, and \texttt{create\_\allowbreak piecewise\_quadratic\_hfun}
builds a max-of-quadratics surrogate from problem-specific $(Q_j, z_j, c_j)$
triples, as used in~\cite{Larson2020}. For both, we reuse the exact
problem data that accompanied those papers' own benchmark studies (the first
random instance/seed for each of the 53 problems), rather than generating our
own. This yields $53 \times 8 = 424$ benchmark problem instances.

Finally, we comment that many of these $h$ functions are unbounded for many of the $F$ instances on
an unconstrained domain.
For example, $h(z) = \min\{z_1, \ldots, z_p\}$ composed with an unconstrained $F$ has
no finite lower bound. To address this, we impose box constraints
$\ell \le x \le u$ around $x_0$ for all problems. 
We employ a simple, deterministic box centered at $x_0$, with
a per-coordinate radius of
\[
r_i = \max(20, \, 2 |x_{0,i}|),
\]
so $\ell_i = x_{0,i} - r_i$ and $u_i = x_{0,i} + r_i$. 
The $\max$ with $20$
matters because several Mor{\'e}--Wild starting points have components at or
near $0$, and so a purely multiplicative radius would collapse to (almost)
zero; scaling by $|x_{0,i}|$ matters because other starting points have
components in the thousands, where a fixed radius would be far too tight to
let a solver explore the domain. 
The function $F$ is not defined for all values in the above box-constrained
domain. 
We observed overflows in $F$ for Meyer's function (corresponding to problem 17) and
Osborne's second function (problems 36 and 37). For these three problems, we use a
tighter box ($r_i = \max(5, |x_{0,i}|)$) to guarantee that $F$ stays finite, instead of employing the radius we prescribed above.

For the IBCDFO test cases, we created an integration module that wraps the BEAD interface and adapts values, gradient
matrices, and path labels to IBCDFO's existing $h$-function interface.
The
PyGRANSO experiments in the following Section use a similar adapter.

\subsection{Other nonsmooth optimization methods}
\label{sec:pygranso-cases}

To verify that BEAD is useful beyond manifold sampling and the explicit
composition $h(F(x))$, we also coupled its branch-extension gradients to PyGRANSO, a
BFGS--SQP method for general nonsmooth constrained optimization
from the NCVX software package \cite{liang2022ncvx,curtis2017bfgssqp}.  We replaced the single-gradient term in
the local search model by the support function of the BEAD gradient set and used
the corresponding worst-case directional derivative in the weak-Wolfe line
search.  
We retained the original BFGS-SQP method as a baseline and also tested an
intermediate variant that used BEAD gradients only in the stationarity test.

We used two small problems for which all within-tolerance paths can be enumerated.  The
first is the constrained nonsmooth Rosenbrock example
\[
 \min_{x_1,x_2}\;8|x_1^2-x_2|+(1-x_1)^2
 \quad\text{s.t.}\quad \sqrt{2}x_1\le 1,\quad 2x_2\le 1,
\]
initialized at the infeasible kink $(1,1)$.
Both the baseline BFGS-SQP method and the
stationarity-only variant failed before accepting an iterate.
Our extension employing BEAD gradients reached a feasible stationary point in 18 iterations and 42 function
evaluations.  

As a second example, we consider a three-variable one-dimensional total-variation
denoising problem (the fused-lasso signal approximator)
\[
 \min_x\;\sum_{i=1}^3(x_i-b_i)^2+\frac{4}{3}
 \sum_{i=1}^2|x_{i+1}-x_i|,\qquad b=(0,2,0).
\]
Starting from $(1,1,1)$, the baseline required 15 iterations and 67 evaluations to terminate at a point near the solution $(2/3,2/3,2/3)$,
using BEAD only for stationarity required 15 iterations and 43 evaluations, and
our variant of the method required one
iteration and three evaluations. 
This behavior can likely be explained by the structure of the problem;
the initial point is already on the same manifold as the solution to the optimization problem, on which the selection function corresponding to the recorded path is quadratic.
These two problems illustrate how BEAD's path and branch-extension information
ports to a general nonsmooth solver; with only two problems, they are meant to be
illustrative and we do not claim they provide a broad performance comparison between the baseline
and BEAD-based variants.

\section{Results}
We now compare BEAD against hand-coded derivatives to see understand how it
affects optimizer performance and then study its runtime.

\subsection{Hand-coded versus BEAD $h$ functions}
\label{sec:hand-vs-bead}

As a direct consistency check, the automated gamma-example test compares
function values and gradient matrices at generic and tie inputs using relative
tolerance $10^{-6}$ and absolute tolerance $10^{-8}$.  The full benchmark below
provides a broader downstream comparison across all eight pairs of $h$
functions.
The remaining question is whether this equivalence translates to downstream performance: does employing BEAD $h$ functions in a manifold sampling algorithm, instead of the hand-coded ones, change the behavior of the optimization method significantly?
We
ran every $(\text{problem}, h)$ pair from the previous subsection twice --
once with each version of $h$ -- through the same manifold sampling solver,
and compared the results using data profiles~\cite{JJMSMW09}.

A data profile plots, for each solver, the fraction of problems solved
as a function of computational budget, measured in simplex gradients
(i.e., groups of $n+1$ function evaluations, where $n$ is the number of
variables for a given problem).
We consider a problem instance to be solved after $k$ evaluations when the
objective value $f(x_k)$ first satisfies
\[
  f(x_k) \;\le\; f_{\mathrm{best}} + \gate\,(f(x_0) - f_{\mathrm{best}}),
\]
where $f(x_0)$ is the starting value, $f_{\mathrm{best}}$ is the smallest
objective value found by any solver on that instance across all evaluations,
and $\gate > 0$ is the convergence gate: \emph{gate} governs the quality of the final optimizer result, while \emph{tolerance} governs the proximity of input or intermediate values that determines which branches are taken during function evaluation.
The gate is a purely function-value-based criterion, measuring how close each solver comes to the best observed value, and can be computed without gradient
information.

In these experiments we use the local branch-discovery rule with
\texttt{rtol}$=0$, so the reported $\tau$ values are the API parameter
\texttt{atol}, the absolute tolerance controlling which alternative branches are considered.
We investigated three values: $\tau = 10^{-8}$ (our baseline, the tightest),
$10^{-6}$ (looser), and $10^{-4}$ (much looser).
This parameter directly affects BEAD's runtime, because a looser tolerance means BEAD must track
more within-tolerance paths. 
The set of discovered paths also affects the downstream optimizer behavior in a problem-dependent and scale-dependent way.
A good choice of $\tau$ is therefore specific to $h$ and to how $h$ is
implemented; the three values with which we experiment are chosen to illustrate this sensitivity,
not to prescribe a general value of $\tau$ useful for all imaginable $h$ functions.

\begin{figure}[h!]
  \centering
  \includegraphics[width=\linewidth]{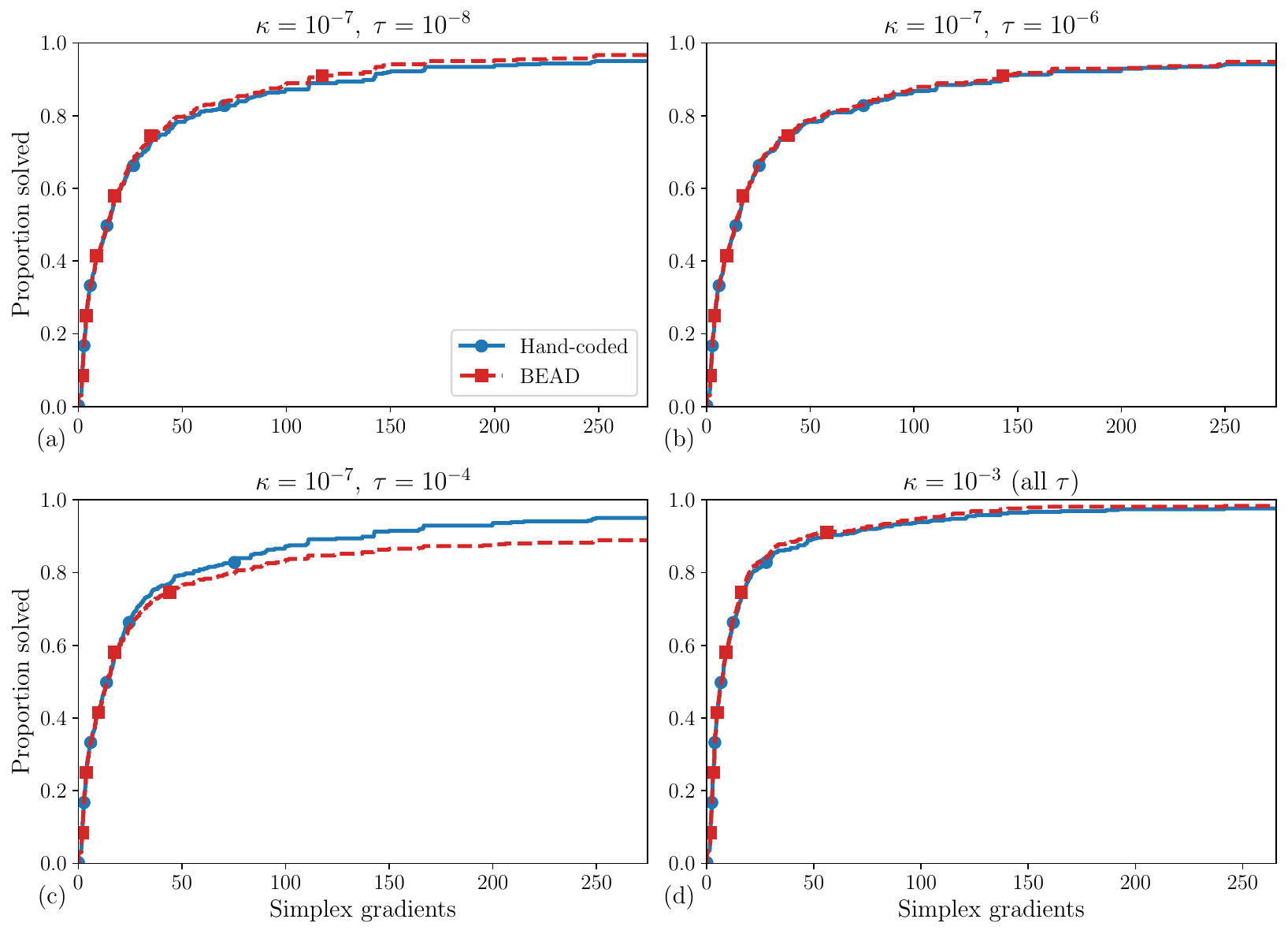}
  \caption{Data profiles comparing hand-coded (solid blue) and BEAD (dashed red)
    $h$ functions across all $424$ (problem, $h$) instances.
    (a)--(c)~Tight convergence gate $\gate = 10^{-7}$ at branch-discovery
    tolerances $\tau = 10^{-8}$, $10^{-6}$, and $10^{-4}$, respectively.
    (d)~Loose gate $\gate = 10^{-3}$: the profiles are indistinguishable
    across all three $\tau$ values, so a single representative is shown.}
  \label{fig:data-profile-all-tols}
\end{figure}

Figure~\ref{fig:data-profile-all-tols} shows data profiles across all three
tolerance values and both convergence gates.
At the loose gate $\gate = 10^{-3}$ (panel~(d)), the two curves are visually
indistinguishable regardless of $\tau$, solving about $93\%$ of the $424$
problem instances within the iteration budget; the three tolerance settings
produce essentially the same profile, so only one representative is shown.
At the tighter gate $\gate = 10^{-7}$ (panels~(a)--(c)), differences become
visible: at $\tau = 10^{-8}$, BEAD and hand-coded performance are comparable; as $\tau$ is loosened, the BEAD curve drops below
hand-coded, reflecting the looser branch-discovery tolerance.
Hand-coded performance is largely insensitive to $\tau$ across this range.

The end result is that replacing a hand-derived $h$ with its BEAD counterpart
does not meaningfully degrade the performance of manifold sampling at the loose
gate or at our tightest tolerance $\tau=10^{-8}$; performance degrades mildly
at the tightest gate as $\tau$ is loosened, but mirrors the coarser
branch-discovery tolerance.
This robustness check supports using automatically generated BEAD path and
branch-extension gradient information in place of hand-coded derivatives.

\subsection{Timing}
\label{sec:timing}

In this section, we evaluate the runtime of BEAD, compared with the hand-coded gradients.
BEAD's branch discovery incurs bookkeeping overhead in each evaluation, and its
current implementation records branches on every call rather than compiling
it once with JAX's JIT.

Table~\ref{tab:timing-hand-vs-jax-all-tols} reports the median wall-clock time per function evaluation for
each $h$, across three values of the branch-discovery tolerance $\tau$.
Both hand-coded and BEAD derivatives are run at the same matched $\tau$ in each column group,
and all timings are restricted to runs that terminated normally (ran the full evaluation budget)
for both versions, so the comparison is over the same amount of work.

At $\tau=10^{-8}$, the BEAD derivative is slower than the hand-coded derivative by a median factor of about $3\times$.
As $\tau$ is loosened from $10^{-8}$ to $10^{-6}$ to $10^{-4}$,
the median ratio increases ($3.1\times \to 3.5\times \to 4.6\times$).
Both implementations may return more gradients at looser tolerances, but their
work differs: the hand-coded routines construct gradient columns directly,
whereas BEAD records additional paths and replays each branch extension through
JAX AD.  The table shows that hand-coded median times remain comparatively
stable while BEAD times increase.
In the simulation-based optimization settings this benchmark represents
(\Cref{sec:test-cases}), $F$ is far more expensive to evaluate than
any computation of derivatives of paths through $h$ and so this overhead is negligible.

\begin{table}[htbp]
  \centering
  \caption{Median wall-clock time per function evaluation (ms), hand-coded vs.\ BEAD $h$ functions, across all dfo rows with normal termination, at matched branch-discovery tolerance $\tau$.}
  \label{tab:timing-hand-vs-jax-all-tols}
  \resizebox{\linewidth}{!}{%
  \begin{tabular}{l@{\hskip 6pt}rr@{\hskip 4pt}r@{\hskip 8pt}rr@{\hskip 4pt}r@{\hskip 8pt}rr@{\hskip 4pt}r}
    \toprule
    & \multicolumn{3}{c}{$\tau = 10^{-8}$} & \multicolumn{3}{c}{$\tau = 10^{-6}$} & \multicolumn{3}{c}{$\tau = 10^{-4}$} \\
    \cmidrule(lr){2-4} \cmidrule(lr){5-7} \cmidrule(lr){8-10}
    $h$ function & hand & BEAD & ratio & hand & BEAD & ratio & hand & BEAD & ratio \\
    \midrule
    censored\_L1\_loss & 9.66 & 35.74 & 3.70 & 7.31 & 40.86 & 5.59 & 14.60 & 35.97 & 2.46 \\
    piecewise\_quadratic & 7.94 & 25.16 & 3.17 & 13.16 & 24.53 & 1.86 & 5.84 & 33.10 & 5.67 \\
    max\_plus\_qvp & 7.74 & 30.39 & 3.93 & 7.81 & 34.12 & 4.37 & 7.17 & 30.92 & 4.31 \\
    one\_norm & 11.26 & 30.07 & 2.67 & 13.75 & 38.22 & 2.78 & 15.74 & 45.24 & 2.87 \\
    pw\_maximum & 7.02 & 23.29 & 3.31 & 7.72 & 25.10 & 3.25 & 7.08 & 30.08 & 4.25 \\
    pw\_maximum\_sq & 8.56 & 24.00 & 2.80 & 8.65 & 29.11 & 3.36 & 6.97 & 36.15 & 5.18 \\
    pw\_minimum & 4.29 & 15.29 & 3.56 & 6.29 & 20.88 & 3.32 & 3.00 & 25.56 & 8.53 \\
    pw\_minimum\_sq & 5.07 & 18.11 & 3.57 & 5.55 & 22.61 & 4.07 & 4.05 & 26.06 & 6.43 \\
    \midrule
    overall (median) & 7.84 & 24.58 & 3.14 & 7.76 & 27.11 & 3.49 & 7.03 & 32.01 & 4.56 \\
    \bottomrule
  \end{tabular}}
\end{table}

\section{Conclusions}

We introduced BEAD, a Python package that collects gradient information for a set of paths through programs describing piecewise-differentiable functions evaluated at or near kinks. 
The package uses JAX, a popular Python framework, for the differentiation of each individual path.

We demonstrate the usefulness of BEAD on two benchmarks: a manifold-sampling benchmark where hand-coded manifold-identification logic was replaced with BEAD, resulting in comparable solver performance for objectives of the form $h(F(x))$; and PyGRANSO examples, showing that the same path and branch-extension gradient information can support a general nonsmooth
optimizer.

The set BEAD returns is not in general a classical generalized derivative;
\Cref{prop:tau-zero} gives conditions on the program and the enumeration policy
under which it contains the Bouligand subdifferential.

The main limitation of BEAD (as well as its hand-coded alternatives) is path-set growth. If a program contains multiple branch points with nearby alternatives (for example, \texttt{max(a, b)} with $a\approx b$), the total number of paths grows exponentially with the number of such branch points. We present one strategy specific to the \texttt{abs} function to mitigate this path explosion, but future work is needed to develop similar policies
for other primitives, and evaluate them in terms of
computational cost and the gradient information that they preserve.

\section*{Acknowledgments}
This work was supported in part by the U.S.~Department of Energy, Office of
Science, Office of Advanced Scientific Computing Research, Scientific Discovery
through Advanced Computing (SciDAC) Program through the FASTMath Institute
under Contract No. DE-AC02-06CH11357.

The initial version of the BEAD software prototype was written by hand, but subsequently refactored, and its API and internal data structures redesigned with the help of AI coding agents. The test cases and experiments were selected by the authors, and testing and plotting scripts developed with the help of AI agents. AI agents were also used in the preparation of early drafts of the user manual and the manuscript, and subsequently heavily revised manually by the authors. The authors take full responsibility for the contents of the manuscript and software submission.

\paragraph{Software availability.}
The BEAD source code, tests, and examples are available in the project
repository~\cite{bead} under the MIT License.
\printbibliography

\end{document}

%% file: figures/diagram_styles.tex
\newlength{\beadgap}
\newlength{\beadpadding}
\newlength{\beadpanelwidth}
\definecolor{beadkeyword}{RGB}{45,86,155}
\definecolor{beadmodule}{RGB}{0,115,120}
\definecolor{beadfunction}{RGB}{110,55,145}
\tikzset{
  bead diagram/.style={font=\scriptsize, node distance=\beadgap,
    every node/.style={outer sep=0pt}},
  bead panel/.style={draw=black!25, rounded corners=3pt, fill=black!1,
    inner sep=\beadpadding},
  bead title/.style={font=\scriptsize\bfseries, text=black!80,
    text width=\beadpanelwidth, align=flush center, inner sep=0pt,
    execute at begin node=\strut, execute at end node=\strut},
  bead box/.style={draw=black!55, rounded corners=2pt, fill=white,
    align=flush center, inner xsep=2pt, inner ysep=3pt, minimum height=7mm},
  bead wide/.style={bead box, text width=\dimexpr\beadpanelwidth-4pt-\pgflinewidth\relax},
  bead code/.style={bead wide, fill=black!3, align=flush left, font=\ttfamily\scriptsize},
  bead decision/.style={diamond, draw=blue!60!black, fill=blue!9,
    aspect=2, align=flush center, inner sep=1pt, font=\scriptsize\bfseries},
  bead data/.style={fill=orange!10, draw=orange!65!black},
  bead replay/.style={fill=purple!8, draw=purple!60!black},
  bead jax/.style={fill=teal!9, draw=teal!60!black},
  bead output/.style={fill=green!8, draw=green!45!black},
  bead note/.style={font=\scriptsize, text=black!65, align=flush center,
    text width=\beadpanelwidth, inner sep=0pt},
  bead connector label/.style={bead note, align=flush left,
    text width=\dimexpr\beadpanelwidth/2-\beadgap\relax},
  bead flow/.style={-{Latex[length=1.5mm]}, semithick, draw=black!65},
  bead default/.style={bead flow, draw=blue!70!black},
  bead nearby/.style={bead flow, dashed, draw=orange!80!black},
  bead selected/.style={bead flow, thick, draw=purple!75!black},
  bead grid/.style={ampersand replacement=\&, inner sep=0pt, outer sep=0pt,
    row sep={\beadgap,between borders},
    column sep={\beadgap,between borders}},
  bead badge/.style={bead wide, minimum height=0pt, font=\scriptsize\bfseries},
  bead include/.style={bead badge, draw=green!55!black, fill=green!10},
  bead exclude/.style={bead badge, draw=red!65!black, fill=red!8}}

%% file: figures/record_replay_design.tex
\begingroup
\setlength{\beadpanelwidth}{\dimexpr(\linewidth-6\beadpadding-2\beadgap-\pgflinewidth)/3\relax}
\newcommand{\programPanel}{%
\begin{tikzpicture}[bead diagram]
  \node[bead title] (title) {(a) User program\\and values at $x_0$};
  \node[bead code, below=of title] (code) {
    \textcolor{beadkeyword}{def} $f(x)$:\\
    \quad\textcolor{beadkeyword}{return}
    (\textcolor{beadmodule}{bnp}.\textcolor{beadfunction}{maximum}(\\
    \qquad $g_1(x),g_2(x)$)\\
    \qquad + \textcolor{beadmodule}{bnp}.\textcolor{beadfunction}{abs}($h(x)$))};
  \matrix[bead grid, below=of code, column sep={.5\beadgap,between borders}] (values) {
    \node[bead box, fill=black!3] (g1) {$g_1(x_0)$\\$=1.00$}; \&
    \node[bead box, fill=black!3] (g2) {$g_2(x_0)$\\$=0.97$}; \&
    \node[bead box, fill=black!3] (h) {$h(x_0)$\\$=0.02$}; \\
  };
  \path (values.south) ++(0,-\beadgap) coordinate (choicetop);
  \path (g1.south) -- (g2.south) coordinate[midway] (maxcenter);
  \node[bead decision, anchor=north] (maxchoice)
    at (maxcenter |- choicetop) {maximum};
  \node[bead decision] (absnode) at (h.south |- maxchoice.center) {abs};
  \coordinate (choicebottom) at (maxchoice.south -| values.south);
  \node[bead box, fill=black!2, circle, inner sep=0pt, minimum size=5mm,
    below=of choicebottom] (plus) {$+$};
  \node[bead wide, bead output, below=of plus] (out) {$f(x_0)$};
  \node[bead note, below=of out] (tolerance)
    {$|g_1(x_0)-g_2(x_0)|$\\$=0.03\leq\tau$\\$|h(x_0)|=0.02\leq\tau$};
  \draw[bead default] (g1.south) -- (maxchoice.north west);
  \draw[bead nearby] (g2.south) -- (maxchoice.north east);
  \draw[bead nearby] (h) -- (absnode);
  \draw[bead flow] (maxchoice.south) -- (plus.north west);
  \draw[bead flow] (absnode.south) -- (plus.north east);
  \draw[bead flow] (plus) -- (out);
\end{tikzpicture}}
\newcommand{\recordPanel}{%
\begin{tikzpicture}[bead diagram]
  \matrix[bead grid, column sep={.5\beadgap,between borders},
    nodes={bead box,
      text width=\dimexpr(\beadpanelwidth-1.5\beadgap)/4-4pt-\pgflinewidth\relax}] (paths) {
    \node[bead data] (p1) {$p_1$\\$g_1$\\$+h$}; \&
    \node[bead data] (p2) {$p_2$\\$g_1$\\$-h$}; \&
    \node[bead data] (p3) {$p_3$\\$g_2$\\$+h$}; \&
    \node[bead replay, thick] (p4) {$p_4$\\$g_2$\\$-h$}; \\
  };
  \path (p1.north) -- (p2.north) coordinate[midway] (pair1);
  \path (p3.north) -- (p4.north) coordinate[midway] (pair2);
  \node[bead decision, above=of pair1] (abs1) {abs};
  \node[bead decision, above=of pair2] (abs2) {abs};
  \path (abs1.north) -- (abs2.north) coordinate[midway] (branchrow);
  \node[bead decision, above=of branchrow] (maximum) {maximum};
  \node[bead wide, bead data, above=of maximum] (record) {record at $x_0$};
  \node[bead title, above=of record] (title) {(b) Recording\\and \texttt{PathSet}};
  \node[bead wide, bead data, below=of paths] (factor)
    {Compact representation:\\$\{g_1,g_2\}\times\{+h,-h\}$\\four resolved paths};
  \draw[bead flow] (record) -- (maximum);
  \draw[bead default] (maximum) -- (abs1);
  \draw[bead nearby] (maximum) -- (abs2);
  \draw[bead default] (abs1) -- (p1);
  \draw[bead nearby] (abs1) -- (p2);
  \draw[bead default] (abs2) -- (p3);
  \draw[bead selected] (abs2) -- (p4);
  \draw[bead flow] (paths.south) -- (factor);
\end{tikzpicture}}
\newcommand{\connectorLabel}[3]{%
  \node[bead connector label, anchor=north west] (#2)
    at ([xshift=\beadgap,yshift=-\beadgap]#1.south) {#3};
  \coordinate (#2-bottom) at (#2.south -| #1.south);
}
\newcommand{\replayPanel}{%
\begin{tikzpicture}[bead diagram]
  \node[bead title] (title) {(c) Locked replay\\and ordinary JAX AD};
  \node[bead wide, bead replay, below=of title] (selected)
    {select $p_4=(g_2,-h)$};
  \node[bead code, below=of selected] (locked)
    {$f_{p_4}(x)=g_2(x)-h(x)$\\
     \textcolor{black!55}{\# fixed forward choices}};
  \connectorLabel{locked}{smooth}{smooth branch extension}
  \node[bead wide, bead jax, below=of smooth-bottom] (ad)
    {{\ttfamily\textcolor{beadmodule}{jax}.\textcolor{beadfunction}{value\_and\_grad}}\\
     $(f_{p_4})(x_0)$};
  \node[bead wide, bead output, below=of ad] (out)
    {$\bigl(g_2(x_0)-h(x_0),$\\$\nabla g_2(x_0)-\nabla h(x_0)\bigr)$};
  \connectorLabel{out}{batchnote}{other paths may be batched}
  \node[bead wide, bead jax, below=of batchnote-bottom] (batch)
    {$p_1,\ldots,p_4$ $\longrightarrow$\\
     {\ttfamily\textcolor{beadmodule}{vmap}(\textcolor{beadfunction}{value\_and\_grad})}\\
     $\longrightarrow$ {\ttfamily\textcolor{beadmodule}{jit/cache}}};
  \draw[bead selected] (selected) -- (locked);
  \draw[bead flow] (locked) -- (ad);
  \draw[bead flow] (ad) -- (out);
  \draw[bead flow] (out) -- (batch);
\end{tikzpicture}}
\newsavebox{\programPanelBox}
\newsavebox{\recordPanelBox}
\newsavebox{\replayPanelBox}
\sbox{\programPanelBox}{\programPanel}
\sbox{\recordPanelBox}{\recordPanel}
\sbox{\replayPanelBox}{\replayPanel}
\begin{tikzpicture}[bead diagram]
  \matrix[bead grid, nodes={anchor=north, inner sep=0pt},
    column sep={\dimexpr2\beadpadding+\beadgap\relax,between borders}] (panels) {
    \node (program) {\usebox{\programPanelBox}}; \&
    \node (recording) {\usebox{\recordPanelBox}}; \&
    \node (replay) {\usebox{\replayPanelBox}}; \\
  };
  \begin{scope}[on background layer]
    \node[bead panel, fit=(program.north west)(program.south east)
      (program.south east |- panels.south)] {};
    \node[bead panel, fit=(recording.north west)(recording.south east)
      (recording.south east |- panels.south)] {};
    \node[bead panel, fit=(replay.north west)(replay.south east)
      (replay.south east |- panels.south)] {};
  \end{scope}
  \matrix[bead grid, below=\dimexpr\beadgap+\beadpadding\relax of panels,
    column sep={\beadgap,between borders}] {
    \node {\tikz{\draw[bead default] (0,0) -- (5mm,0);} default choice}; \&
    \node {\tikz{\draw[bead nearby] (0,0) -- (5mm,0);} within tolerance}; \&
    \node {\tikz{\draw[bead selected] (0,0) -- (5mm,0);} selected replay}; \\
  };
\end{tikzpicture}
\endgroup

%% file: figures/local_input_scaled.tex
\begingroup
\setlength{\beadpanelwidth}{\dimexpr(\linewidth-2\beadpadding-2\beadgap-\pgflinewidth)/3\relax}
\newcommand{\inputView}[5]{%
\begin{tikzpicture}[bead diagram]
  \node[bead title] (heading) {Input-space view};
  \node[inner sep=0pt, below=of heading] (plot) {%
    \begin{tikzpicture}[bead diagram, x=\dimexpr\beadpanelwidth-4pt\relax]
      \fill[blue!12] (0.475,-0.12) rectangle (0.725,0.12);
      \node[anchor=east, inner sep=0pt] (xlabel) at (1,0) {$x$};
      \draw[bead flow] (0,0) -- ([xshift=-\beadgap]xlabel.west);
      \draw[blue!65!black, thick] (0.475,-0.14) -- (0.475,0.14)
        (0.725,-0.14) -- (0.725,0.14);
      \draw[black, thick] (0.6,-0.18) -- (0.6,0.18);
      \draw[#5, thick] (#3,-0.18) -- (#3,0.18);
      \node[bead note, anchor=south] at (0.5,0.28)
        {input ball $|\delta x|\leq\tau=0.1$};
    \end{tikzpicture}};
  \node[bead note, below=of plot] (labels)
    {$x_0=0$ (black)\\\textcolor{#5}{switch $x=#4$}};
  \node[bead wide, below=of labels]
    {$s(x)=#1x+#2$\\$s(x_0)=#2,\quad |s'(x_0)|=#1$};
\end{tikzpicture}}
\newcommand{\ruleView}[5]{%
\begin{tikzpicture}[bead diagram]
  \node[bead title] (heading) {#1};
  \node[bead wide, #2, below=of heading] (rule) {\strut#3\strut};
  \node[#4, below=of rule] {\strut#5\strut};
\end{tikzpicture}}
\newcommand{\comparisonRow}[7]{%
\begin{tikzpicture}[bead diagram]
  \matrix[bead grid, nodes={anchor=north, inner sep=0pt}] (cells) {
    \node {\inputView{#1}{#2}{#3}{#4}{#5}}; \&
    \node {\ruleView{Local rule}{bead data}
      {$|s(x_0)|\leq\tau$\\$#2#6 0.10$}{#7}
      {\localVerdict}}; \&
    \node {\ruleView{Input-scaled rule}{fill=blue!8,draw=blue!60!black}
      {$|s(x_0)|\leq\tau|s'(x_0)|$\\$#2\scaledRelation\scaledBound$}
      {\scaledBadge}{\scaledVerdict}}; \\
  };
\end{tikzpicture}}
\begin{tikzpicture}[bead diagram]
  \node[bead title, text width=\dimexpr\linewidth-2\beadpadding-\pgflinewidth\relax] (ctitle)
    {(a) Contraction: local discovery admits an unreachable choice};
  \def\localVerdict{include both choices}
  \def\scaledVerdict{keep executed choice only}
  \def\scaledRelation{\not\leq}\def\scaledBound{0.02}
  \def\scaledBadge{bead exclude}
  \node[inner sep=0pt, below=of ctitle] (contraction)
    {\comparisonRow{0.2}{0.08}{0.1}{-0.4}{red!70!black}{\leq}{bead include}};
  \node[bead title, text width=\dimexpr\linewidth-2\beadpadding-\pgflinewidth\relax,
    below=\dimexpr\beadgap+2\beadpadding\relax of contraction] (etitle)
    {(b) Expansion: local discovery misses a reachable choice};
  \def\localVerdict{keep executed choice only}
  \def\scaledVerdict{include both choices}
  \def\scaledRelation{\leq}\def\scaledBound{0.50}
  \def\scaledBadge{bead include}
  \node[inner sep=0pt, below=of etitle] (expansion)
    {\comparisonRow{5}{0.20}{0.55}{-0.04}{green!55!black}{\not\leq}{bead exclude}};
  \begin{scope}[on background layer]
    \node[bead panel, fit=(ctitle)(contraction)] {};
    \node[bead panel, fit=(etitle)(expansion)] {};
  \end{scope}
\end{tikzpicture}
\endgroup